\documentclass{mcom-l}

\newtheorem{theorem}{Theorem}[section]

\theoremstyle{definition}

\theoremstyle{remark}

\numberwithin{equation}{section}
\usepackage{amssymb,amsmath,amsthm, amsfonts}
\usepackage{multirow}
\usepackage{mathtools}
\usepackage{subfig}
\usepackage{algorithm,algpseudocode}
\usepackage{caption}
\usepackage{makecell}

\begin{document}

 \title[Preconditioned Iterative Solves in Model Reduction]{Preconditioned Iterative Solves in Model Reduction of Second Order Linear Dynamical Systems}
%\title{Preconditioned Iterative Solves in Model Reduction of Second Order Linear Dynamical Systems}


%    author one information
\author{Navneet Pratap Singh}
\address{Discipline of Computer Science and Engineering, Indian Institute of Technology Indore, India}
\curraddr{}
\email{phd1301201002@iiti.ac.in}
\thanks{}

%    author two information
\author{Kapil Ahuja}
\address{Discipline of Computer Science and Engineering, Indian Institute of Technology Indore, India}
\curraddr{}
\email{kahuja@iiti.ac.in}
\thanks{This work was supported by DAAD grant number A/14/04422 under the IIT-TU9 exchange of faculty program.}

% author 
\author{Heike Fassbender}
\address{Institut Computational Mathematics, Technische Universit\"at Braunschweig, Germany}
\curraddr{}
\email{h.fassbender@tu-braunschweig.de}
\thanks{}

%    \subjclass is required.


\date{}
\dedicatory{}

%    Abstract is required.
\begin{abstract}
Recently a new algorithm for model reduction of second order linear dynamical systems with proportional damping, the Adaptive Iterative Rational Global Arnoldi (AIRGA) algorithm \cite{Bonin20161}, has  been proposed. The  main computational cost of the AIRGA algorithm is in  solving a sequence of linear systems. These linear systems do change only slightly from one iteration step to the next. Here we focus on efficiently solving these systems by iterative methods and the choice of an appropriate preconditioner. We propose the use of relevant iterative algorithm  and the Sparse Approximate Inverse (SPAI) preconditioner. A technique to cheaply update the SPAI preconditioner in each iteration step of the model order reduction process  is given. Moreover, it is shown that under certain conditions the AIRGA algorithm is stable with respect to the error introduced by iterative methods. Our theory is illustrated by  experiments. It is demonstrated that SPAI preconditioned Conjugate Gradient (CG) works well for model reduction of a one dimensional beam model with AIRGA algorithm. Moreover, the computation time of preconditioner with update is on an average $\frac{2}{3}$-rd of the computation time of preconditioner without update. With average timings running into hours for very large systems, such savings are substantial.

\end{abstract}
\maketitle
\section{Introduction}
\label{intro}
\par A continuous time-invariant second order linear dynamical system is of the form
\begin{align}\label{eq1}
	\begin{split}
		M \ddot{x}(t)= & - D \dot{x}(t) - K x(t) + F u(t),  \\
		y(t) =  & \ C_{p} x(t) + C_{v} \dot{x}(t),
	\end{split}  
\end{align}        	   	 
where $M,\ D,\ K \in  \mathbb{R}^{n \times n} $ are mass, damping and stiffness matrices, respectively, $F  \in  \mathbb{R}^{n \times m}, C_{p}, C_{v}  \in  \mathbb{R}^{q \times n}$ are constant matrices. 
In (\ref{eq1}), $ x(t)   \in  \mathbb{R}^{n}$ is the state, $u(t)  \in  \mathbb{R}^{m}$ is the input, and  $y(t)  \in  \mathbb{R}^{q}$ is the output. If $m$ and $q$ both are one, then we have a Single-Input Single-Output (SISO) system. Otherwise ($m$ and $q > 1)$ the system is called Multi-Input Multi-Output (MIMO). We assume the case of proportional damping, i.e., $D=\alpha M + \beta K$, where the coefficients $\alpha$ and $\beta$ are chosen based on experimental results~\cite{BeattieG2005,Bonin20161}. 
For our derivations, the system matrices $M, D$ and $K$ need not hold any specific property (e.g., symmetry, positive definiteness etc.).

It is assumed that the order $n$ of the system (\ref{eq1}) is extremely high. The simulation of large dynamical systems can be unmanageable due to high demands on computational resources, which is the main motivation for model reduction. The goal of model reduction is to produce a low dimensional system that has, as best as possible, the same characteristics as the original system but whose simulation requires significantly less computational effort.  
The reduced  system of (\ref{eq1})  is described by 
\begin{align}\label{TD_r}
	\begin{split}
		\hat{M} \ddot{\hat{x}}(t)= & - \hat{D} \dot{\hat{x}}(t) - \hat{K} \hat{x}(t) + \hat{F} u(t),   \\ 
		\hat{y}(t)=  & \ \hat{C}_{p} \hat{x}(t) + \hat{C}_{v} \dot{\hat{x}}(t),
	\end{split} 	
\end{align}
where $\hat{M}, \ \hat{K}, \ \hat{D}   \in  \mathbb{R}^{r \times r}$, $\hat{F}  \in  \mathbb{R}^{r \times m}$, \ $\hat{C}_{p}, \ \hat{C}_{v} \in  \mathbb{R}^{q \times r}$and $r \ll n$. The damping property of the original system needs to be reflected in the reduced system. That is, $\hat{D}= \alpha \hat M + \beta \hat K$ is required, where $\alpha$ and $\beta$ remain unchanged from the original system.

Model reduction can be done in many ways,
see, e.g., \cite{Antoulas05}. We will focus on a projection based method, specifically Galerkin projection~\cite{Bonin20161}. 
For this a matrix $V \in \mathbb{R}^{n \times r}$ with orthonormal columns is chosen 
and the system (\ref{eq1}) is projected
\begin{align}\label{TD_r1}
	\begin{split}
		& V^{T}(M V\ddot{\hat{x}}(t) + D V \dot{\hat{x}}(t)+KV \hat{x}(t) - F 
		u(t))= 0,\\
		& \hat{y}(t)= \ C_{p} V \hat{x}(t)+C_{v}V\dot{\hat{x}}(t).
	\end{split}
\end{align}
Comparing (\ref{TD_r1}) with (\ref{TD_r}) yields
\begin{align}\label{eq:red_sys}
\begin{split}
	& \hat{M}=V^{T}MV,\ \hat{D}=V^{T}DV,\ \hat{K}=V^{T}KV,\ \hat{F}=V^{T}F, \\ & \hat{C}_{p}=C_{p}V  \ \text{and} \ \hat{C}_{v}=C_{v}V.
\end{split}	
\end{align}
The matrix $V$ can be obtained in many ways, see, e.g.,~\cite{Antoulas05}. 
The focus in this paper will be on the Adaptive Iterative Rational Global Arnoldi (AIRGA) algorithm 
\cite{Bonin20161} in which $V$ is generated by an Arnoldi based approach. 
\par The main contributions of this paper are as follows: Section \ref{sec:1} summarizes the AIRGA model reduction process
which uses a direct solver for solving the linear systems arising in each iteration step.
In Section \ref{sec:Itr_prec}, we discuss the use of  iterative solvers and preconditioners for these  linear systems. 
Preconditioned iterative solvers are a good choice here since they scale well. They have time complexity $\mathcal{O}(n \cdot nnz)$, where $n$ is the size of the system and $nnz$ is the number of nonzeros in system matrices  as compared to  $\mathcal{O}(n^{3})$ for direct solvers~\cite{Saad2003}. The choice of iterative algorithm is problem dependent.  
We show that Sparse Approximate Inverse (SPAI)  preconditioners are well suited for solving the linear systems arising in the model reduction process. These linear systems change at each model reduction iteration, but this change is small. Exploiting this fact we  propose a cheap  preconditioner update. Using an iterative solver introduces additional errors in the computation since the linear systems are not solved exactly. Hence, we discuss the stability of AIRGA in Section~\ref{stablity}.
In Section \ref{Sec:Numerical_exp}, an numerical experiment is given to support our preconditioned iterative solver theory. 
The cheap updates to the SPAI preconditioner, with CG as the underlying iterative algorithm, leads to about $\frac{1}{3}$-rd savings in time. Finally, we give some conclusions and point out future work in Section \ref{sec:conl_fut}.

For the rest of this paper,
$||\cdot||_{F}$ denotes the Frobenius norm, $||\cdot||$  the 2-norm, $||\cdot||_{H_{2}}$  the $H_{2}$-norm, and  $||\cdot||_{H_{\infty}}$  the $H_{\infty}$-norm~\cite{Antoulas05}. 
%If the type of norm is not written, then it is 2-norm.
Also, $qr$ denotes the $QR$ factorization~\cite{trefethen1997numerical}.

%.........section 2
\section{Arnoldi Based Projection Method}
\label{sec:1}
In this section, we first describe how to obtain V such that the first few moments of the transfer functions of the original and the reduced order transfer function are matched. We then state the AIRGA algorithm~\cite{Bonin20161} based on this approach.
%\subsection{Arnoldi using Moment Matching}
%\label{sec:2}
\par The transfer function of (\ref{eq1}) is given by 
\begin{align*}
H(s)=(C_{p}+sC_{v})(s^{2}M+sD+K)^{-1}F=(C_{p}+sC_{v})X(s),
\end{align*}
where $X(s)=(s^{2}M+sD+K)^{-1}F$ is the state variable in frequency domain. The power series expansion of state variable $X(s)$ around expansion point $s_{0} \in \mathbb{R}$ is given as~\cite{Wang2002}
\begin{align}\label{eq:LT_eq1}
X(s)=\sum\limits_{j=0}^{\infty} X^{(j)}(s_{0})(s-s_{0})^{j},
\end{align}
where,
\begin{align}\label{eq:Moment_orig}
%\left \{ \begin{aligned}
\begin{split}
X^{(0)}(s_{0})=& \ (s_{0}^{2}M+s_{0}D+K)^{-1}F,  \\ 
X^{(1)}(s_{0})=& \ (s_{0}^{2}M+s_{0}D+K)^{-1}(-(2s_{0}M+D))X^{(0)}(s_{0}),  \\ X^{(2)}(s_{0})=& \ (s_{0}^{2}M+s_{0}D+K)^{-1}[-(2s_{0}M+D)X^{(1)}(s_{0})-MX^{(0)}(s_{0})],  \\ 
\vdots 
\\ X^{(j)}(s_{0})= & \ (s_{0}^{2}M+s_{0}D+K)^{-1}[-(2s_{0}M+D)X^{(j-1)}(s_{0})  -MX^{(j-2)}(s_{0})].
%\end{aligned}
%\right.
\end{split}
\end{align}
Here, $X^{(j)}(s_{0})$ is called the $j^{th}$-order system moment of $X(s)$ at $s_{0}$.
\par Similarly, the transfer function of the reduced system (\ref{TD_r}) is given by 
\begin{align*}
\hat{H}(s)=(\hat{C}_{p}+s\hat{C}_{v})\hat{X}(s),
\end{align*}
where $\hat{X}(s)=(s^{2}\hat{M}+s\hat{D}+\hat{K})^{-1}\hat{F}.$
The power series expansion of the reduced state space $\hat{X}(s)$ around expansion point $s_{0} \in \mathbb{R}$ is
\begin{align}\label{eq:LT_r}
\hat{X}(s)=\sum\limits_{j=0}^{\infty} \hat{X}^{(j)}(s_{0})(s-s_{0})^{j}.
\end{align}
Here, $\hat{X}^{(j)}(s_{0})$ is defined analogoulsy to the $X^{(j)}(s_{0})$. It is called the $j^{th}$-order system moment of $\hat{X}(s)$ at $s_{0}$.
\par The goal of moment-matching approaches is to find a reduced order model such that the first few moments of (\ref{eq:LT_eq1}) and (\ref{eq:LT_r}) are matched, that is, $X^{(j)}(s_0) = \hat X^{(j)}(s_0)$ for $j = 0, 1, 2, \ldots, t $ for some $t.$
\par Define 
\begin{align*}
	\mathcal{P}_{1} = & - (s_{0}^{2}M+s_{0}D+K)^{-1}(2s_{0}M+ D), \\
	\mathcal{P}_{2} = &-(s_{0}^{2}M+s_{0}D+K)^{-1}M, \\
	\mathsf{Q} = & \ (s_{0}^{2}M+s_{0}D+K)^{-1}F,
\end{align*}
then from (\ref{eq:Moment_orig}) we have 
\begin{align*}
	X^{(0)}(s_{0}) = & \  \mathsf{Q},  \\
	X^{(1)}(s_{0}) = & \ \mathcal{P}_{1}X^{(0)}(s_{0}), \qquad \mathrm{and}\\
	X^{(j)}(s_{0})= & \ \mathcal{P}_{1}  X^{(j-1)}(s_{0}) +  \mathcal{P}_{2}X^{(j-2)}(s_{0})
\end{align*}
for $j \ge 2.$

The second order Krylov subspace \cite{BaiS05} is defined as
\[\mathbb{G}_{j}(\mathcal{P}_{1},\ \mathcal{P}_{2},\ \mathsf{Q}) = \text{span} \{\mathsf{Q},\ \mathcal{P}_1 \mathsf{Q}, 
\ (\mathcal{P}_1^2+  \mathcal{P}_2)  \mathsf{Q}, \ \ldots, \\  \mathcal{S}_j(\mathcal{P}_1,\  \mathcal{P}_2)  \mathsf{Q}\},\] 
where $\mathcal{S}_j( \mathcal{P}_1, \  \mathcal{P}_2) =  \mathcal{P}_1 \cdot \mathcal{S}_{j-1}(\mathcal{P}_1, \ \mathcal{P}_2) +  \mathcal{P}_2 \cdot \mathcal{S}_{j-2}( \mathcal{P}_1,\  \mathcal{P}_2)$ for $j \ge 2.$

Let $\tilde{K}=(s_{0}^{2}M+s_{0}D+K)$. For the special case of proportionally damped second-order systems, it has been observed in \cite{BeattieG2005}
\begin{align*}
\mathbb{G}_{j}( \mathcal{P}_{1},\  \mathcal{P}_{2},\ \mathsf{Q}) 
& =   \mathbb{G}_{j}(-\tilde{K}^{-1}(2s_{0}M+D),\ -\tilde{K}^{-1}M,\ \tilde{K}^{-1}F), \\ 
& = \mathbb{G}_{j}(-\tilde{K}^{-1} ((2s_{0} +\alpha)M+\beta K),\ -\tilde{K}^{-1}M,\ \tilde{K}^{-1}F),\\
& = \mathbb{K}_{j}( \mathcal{P}_1,\  \mathsf{Q}),
\end{align*} 
where $\mathbb{K}_{j}( \mathcal{P}_1,\  \mathsf{Q})$ is the standard block Krylov subspace%. That is,
\begin{align*}
\mathbb{K}_{j}( \mathcal{P}_1, \mathsf{Q})= \text{span} \{ \mathsf{Q}, \  \mathcal{P}_1  \mathsf{Q}, \  \mathcal{P}_1^{2}  \mathsf{Q}, \ \ldots, \  \mathcal{P}_1^{j-1} \mathsf{Q} \}.
\end{align*}
Thus, we need a good basis of $\mathbb{K}_{j}( \mathcal{P}_1,\  \mathsf{Q})$. This can be obtained efficiently by, e.g., the block or the global Arnoldi algorithm~\cite{Saad2003,JbiMS99,Sad93}.

 The AIRGA algorithm, as proposed in~\cite{Bonin20161}, is one of the latest methods based on the global Arnoldi method.
It is given in Algorithm~\ref{Algo:AIRGA}. In this method, moment matching is done at multiple expansion points $s_{i},\ i=\{1,\ \ldots,\ l\},$ rather than just at $s_{0}$ as earlier. This ensures a better reduced model in the entire frequency domain of interest. 

The initial selection and further the computation of expansion points has been discussed in~\cite{Bonin20161} and~\cite{grimme1997phd}. We adopt the choices described in Section 5.0.1 of~\cite{Bonin20161}. The initial expansion points could be either real or imaginary, both of which have their merits. This is problem dependent and discussed in results section. 
After the first AIRGA iteration, the expansion points are chosen from the eigenvalues of the quadratic eigenvalue problem $\lambda^2 \hat{M} + \lambda \hat{D} + \hat{K}$ (at line 33). 

The method is adaptive, i.e., it automatically chooses the number of moments to be matched at each expansion point $s_i.$ This is controlled by the while loop at line 9.
The variable $j$ stores the number of moments matched. The upper bound on $j$ is  $\lceil r_\text{max}/m \rceil$, where $r_\text{max}$ is the maximum dimension to which we want to reduce the state variable (input from the user), and $m$ is the dimension of the input; see \cite{Bonin20161} for a detailed discussion. At exit of this while loop, $J = j$.
\begin{algorithm}[!]
	\normalsize
	\caption{Adaptive Iterative Rational Global Arnoldi Algorithm \cite{Bonin20161}}
	\label{Algo:AIRGA}
	\begin{algorithmic}[1]		
		\\ Input:  \{$M,\ D,\ K,\ F,\ C_{p},\ C_{v}$,\ $r_\text{max}$; \ $S$ is the set initial expansion points $s_{i},\ i=1,\ \ldots, \ l$\}
		%\\ $J = r / size(F,2) $;
		%	\While {$||\hat{H}_{old}-\hat{H}_{new}||_{H_{2}}\leq tol$}
		\State $z = 1$
		\While{no convergence}
		\For  {$ \text{each}\ s_{i} \in S$} \qquad %\qquad \small{\textit{* /Initialize */}}
		%\State $R^{(-1)}(s_{i})=0$
		\State $X^{(0)}(s_{i})= (s_{i}^{2}M+s_{i}D+K)^{-1}F$ 
		\State Compute $QR = qr(X^{(0)}(s_{i}))$,\ 
		$X^{(0)}(s_{i})=Q$ 
		\EndFor
		\State j = 1 
		\While{no convergence and $j \le \lceil r_\text{max} / m \rceil $}
		\State {Let $\sigma_{j}$ be expansion point corresponding to maximum moment error of \Statex \qquad \quad reduced  system} at $s_{i}$
		\State $V_{j}=X^{(j-1)}(\sigma_{j})/||X^{(j-1)}(\sigma_{j})||_{F}$
		\For {$i=1,\ \ldots,\ l$} 
		\If {($s_{i}==\sigma_{j}$)} 
		\State $X^{(j)}(s_{i})=-(s_{i}^{2}M+s_{i}D+K)^{-1}MV_{j}$
		\Else $ \ X^{(j)}(s_{i})=X^{(j-1)}(s_{i})$
		\EndIf
		\For {$t=1,\ 2,\ \ldots,\ j$}
		\State $\gamma_{t,j}(s_{i})=\text{trace}(V_{t}^{H}\cdot X^{(j)}(s_{i}))$ \State $X^{(j)}(s_{i})=X^{(j)}(s_{i})-\gamma_{t,j}(s_{i})V_{t}$ 
		\EndFor
		\EndFor
		\State $W_i =  X^{(j)}(s_{i})/||X^{(j)}(s_{i})||_{F}$ for $i = \{1, \ldots, l\}$.
		\State $\tilde{W}=[W_{1},\ W_{2},\ \ldots,\ W_{l}]$.
		\State Compute $\check{W}Y = qr(\tilde{W}),\ W=\check{W}$ 
		\State  Compute reduced system matrices $\hat{M}, \hat{D},$ and $ \hat{K} $ with $V = W$ as in~(\ref{eq:red_sys})
		\State $\hat{H}_{\text{Int}}=(\hat{C}_{p}+\sigma_{j}\hat{C}_{v})(\sigma_{j}^{2}\hat{M}+\sigma_{j}\hat{D}+\hat{K})^{-1}\hat{F}$
		\State j = j+1
		\EndWhile 
		\State {Set $J = j$ and pick $\sigma_{J}$ corresponding to  maximum moment error of  reduced \Statex \quad \ \  system at $s_{i}$} 
		\State $V_{J}=X^{(J-1)}(\sigma_{J})/||X^{(J-1)}(\sigma_{J})||_{F}$ and $\tilde{V}=[V_{1},\ V_{2},\ \ldots,\ V_{J}]$.
		\State Compute $\check{V}Y = qr(\tilde{V}),\ V=\check{V}$ 
		\State  Compute reduced system matrices $\hat{M}, \hat{D},$ and $ \hat{K} $ with $V$ as in~(\ref{eq:red_sys}), \Statex \quad \ \  and take $M = \hat{M}$,  
			$D = \hat{D}$,  $K = \hat{K}$, for the next iteration.
		\State  Choose new expansion points $s_{i}$
		\State $\hat{H}=(\hat{C}_{p}+\sigma_{J}\hat{C}_{v})(\sigma_{J}^{2}\hat{M}+\sigma_{J}\hat{D}+\hat{K})^{-1}\hat{F}$
		\State $z = z + 1$
		\EndWhile
		\\ Compute the remaining reduced system matrices $\hat{F}, \ \hat{C}_{p}, \ \hat{C}_{v}$ with $V$ as in~(\ref{eq:red_sys})
	\end{algorithmic}
\end{algorithm}

At line 3, no convergence implies that the $H_2$ norm of the difference between two consecutive reduced systems, computed at line 34, is greater than a certain tolerance. Similarly, at line 9, no convergence implies that the $H_2$ norm of the difference between two consecutive intermediate reduced systems, computed at line 26, is greater than a certain tolerance.     

This algorithm requires solving a linear system at line 5 and 14. As the $s_i$ change in each iteration step, the linear systems to be solved change in each iteration step. As discussed in Section~\ref{intro}, since solving such systems by direct methods is quite expensive, we propose to use iterative methods. As the change in the $s_i$ will be small (at least after the first iteration step), we can develop a cheap update of the necessary preconditioner. 

\section{Preconditioned Iterative Method}
\label{sec:Itr_prec}

\par There are two types of methods for solving linear systems of equations; a) direct methods and b) iterative methods. For  large systems, direct methods are
not preferred because they are too expensive in terms of storage and operation. On the other hand, iterative methods require less storage and operations than direct methods. 
For a large linear system $Ax=b,$
with $A  \in \mathbb{R}^{n \times n}$ and $b  \in  \mathbb{R}^{n}$,  an iterative method finds a sequence of solution vectors $x_{0},\ x_{1},\ \ldots,$ $\ x_{k}$ which (hopefully) converges to the desired solution. Krylov subspace based methods are an important and popular class of iterative methods. If $x_{0}$ is the initial solution and $r_{0}=b-Ax_{0}$ is the initial residual, then Krylov subspace methods find the approximate solution by projecting onto the Krylov subspace 
\begin{align*}
\mathbb{K}_{k}(A,\ r_{0})=\text{span}\{r_{0}, \ Ar_{0},\ A^{2}r_{0},\ \ldots,\ A^{k-1}r_{0}  \}.
\end{align*}
\par There are many types of Krylov subspace algorithms~\cite{Saad2003}. Some popular ones include Conjugate Gradient (CG), Generalized Minimal Residual (GMRES), Minimum Residual (MINRES), and BiConjugate Gradient (BiCG). Block versions of these algorithms do exist. The choice of algorithm is problem dependent. In results section (Section~\ref{Sec:Numerical_exp}), the coefficient matrices arising from the problem are symmetric positive definite~(SPD). Hence, we use CG algorithm, which is ideal for such systems.       

In Krylov subspace methods, the conditioning of the system is very important. ``Conditioning pertains to the perturbation behavior of a mathematical problem~\cite{trefethen1997numerical}''. For example, in a well-conditioned problem, a small perturbation of the input leads to a small change in the output. This is not guaranteed for an ill-conditioned problem, where a small perturbation in the input may change the output drastically~\cite{trefethen1997numerical}. Preconditioning is a technique to well-condition an ill-conditioned problem. We discuss that next. 

Preconditioning is a technique for improving the performance of iterative methods. It transforms a difficult system (ill-conditioned system) to another system with more favorable properties for iterative methods. For example, a preconditioned matrix may have eigenvalues clustered around one. This means that the preconditioned matrix is close to the identity matrix, and hence, the iterative method will converge faster.  For a symmetric positive definite (SPD) system, the convergence rate of  iterative methods depends on the distribution of the eigenvalues of the coefficient matrix. However, for a non-symmetric system, the convergence rate may depend on pseudo-spectra as well ~\cite{trefethen1991pseudospectra,nachtigal1992fast}. 
\par If $M$ is a nonsingular matrix which approximates A; that is, $M \approx A^{-1},$ then the system 
\begin{align*}%\label{Eq:preAM}
MAx=Mb 
\end{align*}
may be faster to solve than the original one. The above system represents preconditioning from left. Similarly, right and split preconditioning is given by two equations below, respectively. 
\begin{align*}
\begin{split}
AM\tilde{x}=b, \ x=M\tilde{x} \ \ \text{and} \  \ 
M_{1}AM_{2}\tilde{x}=M_{1}b, \ x=M_{2}\tilde{x}.
\end{split}
\end{align*}   
The type of preconditioning technique to be used depends on the problem properties as well as on the choice of the iterative solver. For example, for SPD systems,$\ MA, \ AM, \ M_{1}AM_{2}$ all have same eigenvalue spectrum, and hence, left, right and split preconditioners behave the same way, respectively. For a general system, this need not be true~\cite{Benzi2002}.                    
\par Besides making the system easier to solve by an iterative method, a preconditioner should be cheap to construct and apply. Some existing preconditioning techniques include Successive Over Relaxation, Polynomial, Incomplete Factorizations, Sparse Approximate Inverse~(SPAI), and Algebraic Multi-Grid~\cite{Benzi2002, benzi1999comparative, saadchow, Chow1998}. 

%Successive Over Relaxation and Polynomial preconditioners are very basic, and work when the coefficient matrix of the linear system to be solved has nice properties, e.g. it may be diagonally dominance or symmetric.  Our coefficient matrices do not have such properties. Incomplete Factorization based preconditioners are mature, and are commonly used. However, they also work well when the coefficient matrix has some of the above discussed properties, which we do not have here. Algebraic Multi-Grid preconditioners are extremely useful for solving linear systems arising from full discretization of a Partial Differential Equation (PDE). The linear systems here do not satisfy this. 

%Fortunately, SPAI preconditioners are useful here. These preconditioners are known to work well in the most general setting. 
We use SPAI preconditioner here since these~(along with incomplete factorizations) are known to work in the most general setting. Also, SPAI preconditioners are easily parallelizable, hence, have an edge over incomplete factorization based preconditioners~\cite{van2003iterative}.    

In Section~\ref{Sec:SPAI} we summarize the SPAI preconditioner from~\cite{Chow1998} and we discuss the use of SPAI in the AIRGA algorithm. Since the change in the coefficient matrix of the linear system to be solved is small from one step of AIRGA  to the next, we update the preconditioner from one step to the next. This aspect is covered in Section~\ref{Sec:SPAI_U}.     
% end section 3
% subsectio 3
\subsection{Sparse Approximate Inverse (SPAI) Preconditioner}\label{Sec:SPAI}
In constructing a preconditioner $P_{i}^{(z)}$ for a coefficient matrix $\mathcal{K}_{i}^{(z)}$ at the $z^{\text{th}}$ outer AIRGA iteration (i.e., $\mathcal{K}_{i}^{(z)}= s_{i}^{2}M+s_{i}D+K$), we would like $P_{i}^{(z)}\mathcal{K}_{i}^{(z)} \approx I$ (for left preconditioning) and
$\mathcal{K}_{i}^{(z)}P_{i}^{(z)} \approx I$ (for right preconditioning). SPAI preconditioners find $P_{i}^{(z)}$ by minimizing the associated error norm $\left|\left|I-P_{i}^{(z)}\mathcal{K}_{i}^{(z)}\right|\right|$ or $\left|\left|I-\mathcal{K}_{i}^{(z)}P_{i}^{(z)}\right|\right|$ for a given sparsity pattern. If the norm used is Frobenius norm, then the minimization function will be
\begin{align*}
\min_{P_{i}^{(z)} \in S} \left|\left|I-\mathcal{K}_{i}^{(z)}P_{i}^{(z)}\right|\right|_{F},
\end{align*}
where $S$ is a set of certain sparse matrices. The above approach produces a right approximate inverse. Similarly, a left approximate inverse  can be computed by solving the minimization problem $\left|\left|I-P_{i}^{(z)}\mathcal{K}_{i}^{(z)}\right|\right|_{F}$. For non-symmetric matrices, the distinction between left and right approximate inverses is important. There are some situations where it can be difficult to find a right approximate inverse but finding a left approximate inverse can be easy. Whether left or right preconditioning should be used is problem dependent~\cite{Saad2003}. Since the SPAI preconditioner was originally proposed for right preconditioning~\cite{Chow1998}, we focus on the same here. Similar derivation can be done for the left preconditioning as well.
\begin{algorithm}[]
	%\large
	\normalsize
	\caption{: Sparse Approximate Inverse (SPAI) Preconditioner~\cite{Chow1998}}
	\label{ALGO:SPAI}
	\begin{algorithmic}[1]
		\\ Input: \{$\mathcal{K}_{i}^{(z)}, tol$\}
		% \\ Input: \{ $\mathcal{K}_{i}^{(z)} ,n_i$, \textit{lfil} \}
		\\$P_{i}^{(z)} = \alpha I$  \ where \ $\alpha=\dfrac{\text{trace}\left(\mathcal{K}_{i}^{(z)}\right)}{\text{trace}\left(\mathcal{K}_{i}^{(z)}\left(\mathcal{K}_{i}^{(z)}\right)^{T}\right)}$ \ and $n$ is the dimension of $\mathcal{K}_{i}^{(z)}$
%		\For {$outer = 1,\ \ldots,\ n_{0}$}
		%\While {$||I - \mathcal{K}_i P_i||_{F} < 1$}
		\For {$ j = 1,\ \ldots,\ n$}
		\State Define $p_{i}^{(j)}=P_{i}^{(z)}e^{(j)}$
		\State $r = e^{(j)}-\mathcal{K}_{i}^{(z)}p^{(j)}_{i}$
        	\While {$||r|| > tol$}
		% \For {$inner =1, \ \ldots, \ n_{i}$}
        %\While { $||r||$ $<$ 1}
		%\State $r = e^{(j)}-\mathcal{K}_{i}p^{(j)}_{i}$
		\State $d =P_{i}^{(z)}r$
		% \State $t =P_{i}^{(z)}r$
		% \State Choose $d$ to be $t$ with the same pattern as $p_i^{(j)};$
		% \Statex \qquad  \quad If $(nnz(p_i^{(j)})) <$ \textit{lfil}, then add one entry which is the \Statex \qquad \quad largest remaining  entry in absolute value 
		\State $w =\mathcal{K}_{i}^{(z)}d$
		\State $\alpha  = \frac{(r,\ w)}{(w,\ w)}$
		\State $p^{(j)}_{i} =p^{(j)}_{i}+\alpha d$
		\State $r=r-\alpha w$ 	
		\EndWhile
		%\State Update $j^{th}$ column of $P_{i}$ with $p^{(j)}_{i}$
		\EndFor
		%\EndFor	
	\end{algorithmic}
\end{algorithm} 
The minimization problem can be rewritten as 
\begin{align}\label{eq:Lq}
\min\left|\left|I-\mathcal{K}_{i}^{(z)}P_{i}^{(z)}\right|\right|_{F}^{2}= \min\sum_{j=1}^{n} \left|\left|e^{(j)}-\mathcal{K}_{i}^{(z)}p^{(j)}_{i}\right|\right|^{2}_2,
\end{align}
where $p^{(j)}_{i}$ and $e^{(j)}$ are  $j^{th}$ columns of  the $P_{i}^{(z)}$ matrix and \textit{I} (identity matrix), respectively. The minimization problem \eqref{eq:Lq} is essentially just one least squares problem, to be solved for $n$ different right-hand sides. 
Here it is solved iteratively.
\par The algorithm for computing a SPAI preconditioner for right preconditioning is given in Algorithm~\ref{ALGO:SPAI}. The inputs to this algorithm is $\mathcal{K}_{i}^{(z)}$ (coefficient matrix) and $tol$ (stopping residual of the minimization problem for each column);  $tol$ is picked based on experience. This is ALGORITHM 2.5 of~\cite{Chow1998} with two minor differences. 

First, we do not list the code related to sparsity pattern matching (for obtaining a sparse preconditioner) because the goal here is to motivate SPAI update, and for our problems the original matrix is very sparse so the preconditioner stays sparse any ways. This aspect can be easily incorporated. Second, we use a $While$ loop at line 6 of  Algorithm~\ref{ALGO:SPAI} instead of a $For$ loop. This is because with a $For$ loop one has to decide the stopping count in advance (which is chosen heuristically). We use a more certain criteria. That is, residual of the minimization problem for each column less than $tol$. This is linear cost (for each column) and we are doing such computation anyways.  
 
The initial guess for approximate inverse $P_{i}^{(z)}$ is usually taken as $\alpha \textit{I}$ where \\ $\alpha = {\text{trace}\left(\mathcal{K}_{i}^{(z)}\right)}/{\text{trace}\left(\mathcal{K}_{i}^{(z)}\left(\mathcal{K}_{i}^{(z)}\right)^{T}\right)}$ (see line 2). This initial scaling factor $\alpha$ is 
minimizes the spectral radius of $(I-\alpha \mathcal{K}_{i}^{(z)})$~\cite{benzi1999comparative,Chow1998,CarlDMeyer}. 

\par The AIRGA algorithm with the SPAI preconditioner is given in Algorithm~\ref{ALGO:AIRGASPAI}. Here, we only show those parts of AIRGA algorithm that require changes.  
\begin{algorithm}[]
	%	\large
	\normalsize
	\caption{: AIRGA Algorithm with SPAI Preconditioner}
	\label{ALGO:AIRGASPAI}
	\begin{algorithmic}[1]
		%\While {$||\hat{H}_{old}-\hat{H}_{new}||_{H_{2}} \ \leq \ tol$}
		\While{no convergence}
		\For {$i = 1,\ \ldots,\ l $}  
		\State Let $\mathcal{K}_{i}^{(z)}=(s^{2}_{i}M+s_{i}D+K)$ 
		\State Compute preconditioner $P_{i}^{(z)}$ by solving  
		$\min\left|\left|I-\mathcal{K}_{i}^{(z)}P_{i}^{(z)}\right|\right|_{F}^{2}$
		\State Solve $\mathcal{K}_{i}^{(z)}P_{i}^{(z)}\tilde{X}^{(0)}(s_{i})=F$ with ${X}^{(0)}(s_{i})=P_{i}^{(z)}\tilde{X}^{(0)}(s_{i})$ 		 
		\EndFor
		%\While {$||\hat{H}_{j}-\hat{H}_{j+1}||_{H_{2}} \leq tol $}
		\State j = 1
		\While{no convergence and $j \le \lceil r_\text{max} / m \rceil $}
		\For {$i = 1,\ \ldots, \ l$}
		\State Only right hand sides are changing, so above preconditioner 
		$P_{i}^{(z)}$ can 
		\Statex \qquad \qquad \ \ be  applied as it is, i.e.,   
		\Statex \qquad \qquad \ \ 
		Solve $\mathcal{K}_{i}^{(z)}P_{i}^{(z)}\tilde{X}^{(j)}(s_{i})=MV_{j}$ with ${X}^{(j)}(s_{i})=P_{i}^{(z)}\tilde{X}^{(j)}(s_{i})$ 	
		\EndFor
		\EndWhile
		\State j = j+1
		\EndWhile
	\end{algorithmic}
\end{algorithm}

\subsection{SPAI Update Preconditioner}\label{Sec:SPAI_U}
Let $\mathcal{K}_{{old}}=s_{old}^{2}M+s_{old}D+K$ and $\mathcal{K}_{new}=s_{new}^{2}M+s_{new}D+K$ be two coefficient matrices for different expansion points $s_{old}$ and $s_{new}$, respectively. These expansion points can be at the same or different AIRGA iteration. If the difference between $\mathcal{K}_{old}$  and $\mathcal{K}_{new}$ is small, then one can exploit this while building preconditioners for this sequence of matrices. This has been considered in the quantum Monte Carlo setting~\cite{ahuja2011improved} and for model reduction of first 
order linear dynamical systems~\cite{grim2015reusing,Wyatt2012}.  
\par Let $P_{old}$ be a good initial  preconditioner for $\mathcal{K}_{old}$. 
As will be seen, a cheap preconditioner update can be obtained  by asking for $\mathcal{K}_{old}P_{old} \approx \mathcal{K}_{new}P_{new}$, where $old,new =\{1,\ldots,l\}$ and, as earlier, $l$ denotes the number of expansion points.    
Expressing $\mathcal{K}_{new}$ in terms of $\mathcal{K}_{old}$, we get
\begin{align*}
\mathcal{K}_{new}=\mathcal{K}_{old}(I+(s^{2}_{new}-s^{2}_{old})\mathcal{K}_{old}^{-1}M+(s_{new}-s_{old})\mathcal{K}_{old}^{-1}D).
\end{align*}
Now we enforce $\mathcal{K}_{old}P_{old}=\mathcal{K}_{new}P_{new}$ or
\begin{align*}
\begin{split}
\mathcal{K}_{old}P_{old}=& \ \mathcal{K}_{old}(I+(s^{2}_{new}-s^{2}_{old})\mathcal{K}_{old}^{-1}M+(s_{new}-s_{old})\mathcal{K}_{old}^{-1}D) \\ 
& \ \cdot (I + (s^{2}_{new}-s^{2}_{old})\mathcal{K}_{old}^{-1}M + (s_{new}-s_{old})\mathcal{K}_{old}^{-1}D)^{-1}P_{old} \\
= &\ \mathcal{K}_{new}P_{new},
\end{split}
\end{align*}
where     $P_{new}=(I+(s^{2}_{new}-s^{2}_{old})\mathcal{K}_{old}^{-1}M+(s_{new}-s_{old})\mathcal{K}_{old}^{-1}D)^{-1}P_{old}$. 

Let $Q_{new}\approx (I+(s^{2}_{new}-s^{2}_{old})\mathcal{K}_{old}^{-1}M+(s_{new}-s_{old})\mathcal{K}_{old}^{-1}D)^{-1}$,  then the above implies $\mathcal{K}_{old}P_{old}\approx \mathcal{K}_{new}Q_{new}P_{old}$ or $\mathcal{K}_{old}\approx \mathcal{K}_{new}Q_{new}$. This
leads us to the following idea: instead of solving for $P_{new}$ from $\mathcal{K}_{old}P_{old} = \mathcal{K}_{new} P_{new}$, we solve a simpler problem
\begin{equation*}
\min\left|\left|\mathcal{K}_{old}-\mathcal{K}_{new}Q_{new}\right|\right|^{2}_{F}=\min \sum_{j=1}^{n}\left|\left|k^{(j)}_{old}-\mathcal{K}_{new}q^{(j)}_{new}\right|\right|^{2}_2,
\end{equation*}
where $k_{old}^{(j)}$ and $q_{new}^{(j)}$ denote the $j^{th}$ columns of $\mathcal{K}_{old}$ and $Q_{new}$, respectively. Compare this minimization problem with the one in SPAI (Equation (\ref{eq:Lq}) in Section 3.1). Earlier, we were finding the preconditioner $P_{new}$ for $\mathcal{K}_{new}$ by solving $\min \left|\left|I-\mathcal{K}_{new}P_{new}\right|\right|_{F}^{2}$. Here, we are finding the preconditioner $P_{new}$ (i.e., $P_{new}=Q_{new}P_{old}$) by solving $\min\left|\left|\mathcal{K}_{old}-\mathcal{K}_{new}Q_{new}\right|\right|_{F}^{2}$. The second formulation is much easier to solve since in the first $\mathcal{K}_{new}$ could be very different from $I$, while in the second $\mathcal{K}_{new}$ and $\mathcal{K}_{old}$ are similar (change only in the expansion points).

The SPAI algorithm (Algorithm 2) adapted for finding the preconditioner by minimizing this new expression is given in Algorithm~\ref{ALGO:SPAI_U}. The inputs to this algorithm include $\mathcal{K}_{old}$, $\mathcal{K}_{new}$,  and $tol$ (stopping residual of the minimization problem for each column);  $tol$ is picked based on experience. The initial guess for the approximate inverse $Q_{new}$ is usually taken as $\alpha \textit{I}.$ Similar to before, $\alpha$ is chosen to minimize the spectral radius of $(\mathcal{K}_{old}-\alpha \mathcal{K}_{new}):$  
\begin{align*}
\begin{split}
\frac{\partial}{\partial \alpha}\left|\left|\mathcal{K}_{old}-\alpha \mathcal{K}_{new}\right|\right|_{F}^{2}=0 \quad \text{or} \\
\frac{\partial}{\partial \alpha}\left|\left|\mathcal{K}_{old}-\alpha \mathcal{K}_{new}\right|\right|_{F}^{2}=\frac{\partial}{\partial \alpha}\text{trace}(\mathcal{K}_{old}-\alpha \mathcal{K}_{new})^{T}(\mathcal{K}_{old}-\alpha \mathcal{K}_{new}) =0  \quad \text{or} \\
\frac{\partial}{\partial \alpha} \text{trace} [\mathcal{K}_{old}^{T}\mathcal{K}_{old}- \alpha \mathcal{K}_{old}^{T}\mathcal{K}_{new}-\alpha  \mathcal{K}_{new}^{T}\mathcal{K}_{old} +\alpha^{2}\mathcal{K}_{new}^{T}\mathcal{K}_{new}]=0 \quad \text{or} \\
2\alpha \cdot \text{trace} (\mathcal{K}_{new}^{T}\mathcal{K}_{new})= \text{trace}(\mathcal{K}_{old}^{T}\mathcal{K}_{new}+\mathcal{K}_{new}^{T}\mathcal{K}_{old})  \quad \text{or} \\
\alpha =  \frac{1}{2}\cdot\frac{\text{trace}(\mathcal{K}_{old}^{T}\mathcal{K}_{new}+\mathcal{K}_{new}^{T}\mathcal{K}_{old})}{\text{trace}(\mathcal{K}_{new}^{T}\mathcal{K}_{new})}.
\end{split}
\end{align*} 
\begin{algorithm}[!t]
	%\large
	\normalsize
	\caption{: SPAI Update Preconditioner}
	\label{ALGO:SPAI_U}
	\begin{algorithmic}[1]
		\\ Input: \{$\mathcal{K}_{old},\ \mathcal{K}_{new},\ tol$\}
		% \\ Input: \{ $\mathcal{K}_{old},\ \mathcal{K}_{new}, \ n_i$, \textit{lfil} \}
		\State $Q_{new}=\alpha I$  
		\Statex    where \ $\alpha =  \frac{1}{2}\cdot\dfrac{\text{trace}(\mathcal{K}_{old}^{T}\mathcal{K}_{new}+\mathcal{K}_{new}^{T}\mathcal{K}_{old})}{\text{trace}(\mathcal{K}_{new}^{T}\mathcal{K}_{new})}$ \ and $n$ is the dimension of $\mathcal{K}_{new}$
	%	\For {$outer = 1,\ \ldots,\ n_{0}$}
		%\While {$||\mathcal{K}_{i_{old}} - \mathcal{K}_iQ_i||_{F} < 1$}
		\For {$j = 1,\ \ldots,\ n$}
		\State Define $q^{(j)}=Q_{new}e^{(j)}$
		\State $r = k_{old}^{(j)}-\mathcal{K}_{new}q^{(j)}$
		
		\While {$||r|| > tol$}
		% \For {$inner = 1,\ \ldots,\ n_{i}$}
%		\While {$||r|| < 1$}
			\State $d = Q_{new}r$ 
			% \State $t = Q_{new}r$ 
			% \State Choose $d$ to be $t$ with the same pattern as $q^{(j)};$
		% \Statex \qquad  \quad If $(nnz(q^{(j)})) <$ \textit{lfil}, then add one entry which   \Statex \qquad \quad is the largest remaining  entry in absolute value 		
		\State $w =\mathcal{K}_{new}d$ 
		\State $\alpha = \frac{(r,\ w)}{(w,\ w)}$
		\State $q^{(j)} = q^{(j)}+\alpha d$
		\State $r = r - \alpha w$
		\EndWhile
		%\State  Update $j^{th}$ column of $Q_i$ with $q^{(j)}$
		\EndFor
	%	\EndFor	
	\end{algorithmic}
\end{algorithm}  
\begin{algorithm}[!tbp]
	%\large
	\normalsize
	\caption{: AIRGA with SPAI Update Preconditioner}
	\label{ALGO:AIRGASPAI_Update}
	\begin{algorithmic}[1]
		%\While {$||\hat{H}_{old}-\hat{H}_{new}||_{H_{2}} \ \leq \ tol$}
		\State $z = 1$
		\While{no convergence}
		\For {i = 1 to $l$}
		\If {$z == 1$}
		\State $\mathcal{K}_i^{(1)}=(s^{2}_iM+s_iD+K)$
		\State Compute initial $P_i^{(1)} $ by solving $ \min\left|\left|I-\mathcal{K}_i^{(1)}P_i^{(1)}\right|\right|_{F}^2$
		\State Solve $\mathcal{K}_i^{(1)}P_i^{(1)}\tilde{X}^{(0)}(s_i) = F$ with ${X}^{(0)}(s_{i})=P_{i}^{(1)}\tilde{X}^{(0)}(s_{i})$ 	
		%\State Store $\mathcal{K}_{i_{old}} = \mathcal{K}_i$ and $P_{i_{old}} = P_i$
		\Else
		\State $\mathcal{K}_i^{(z)}=(s^{2}_iM+s_iD+K)$
		\State Compute $Q_i^{(z)} $ by solving $ \min\left|\left|\mathcal{K}_i^{(z-1)}-\mathcal{K}_i^{(z)}Q_i^{(z)}\right|\right|_{F}^2$
		%\State Compute next preconditioner $P_i = Q_i P_{i_{old}}$
		\State Solve $\mathcal{K}_i^{(z)}\left[Q_i^{(z)} Q_i^{(z-1)} \ldots Q_i^{(2)} P_i^{(1)}\right]\tilde{X}^{(0)}(s_i) = F$  \Statex \qquad \qquad \ \
		with ${X}^{(0)}(s_i)=\left[Q_i^{(z)} Q_i^{(z-1)} \ldots Q_i^{(2)} P_i^{(1)}\right] \tilde{X}^{(0)}(s_i)$
		% \State Store $\mathcal{K}_{i_{old}} = \mathcal{K}_i$ and $P_{i_{old}} = P_i$	
		\EndIf
		\EndFor
				\State $j = 1$
				\While{no convergence and $j \le \lceil r_\text{max} / m \rceil $}
				\For {i = 1 to $l$}
				\State Only right hand sides are changing, so above  preconditioners can be \Statex \qquad \qquad \ \ applied as it is, i.e.,
				\Statex \qquad \qquad \ \ Solve $\mathcal{K}_i^{(z)}\left[Q_i^{(z)} Q_i^{(z-1)} \ldots Q_i^{(2)} P_i^{(1)}\right] \tilde{X}^{(j)}(s_i)=MV_j$ \Statex \qquad \qquad \ \ with ${X}^{(j)}(s_i)=\left[Q_i^{(z)} Q_i^{(z-1)} \ldots Q_i^{(2)} P_i^{(1)}\right] \tilde{X}^{(j)}(s_i)$
		\EndFor		
		\EndWhile
			\State $j = j+1$
		\EndWhile
		\State $z = z +1$
	\end{algorithmic}
\end{algorithm}

\begin{table}[!b]
	\centering
	\caption{Change in Expansion Points }
	\label{change_in_expansion_pnt}
	\begin{tabular}{|c|c|c|c|c|c|}
		\hline
		\textbf{\begin{tabular}[c]{@{}c@{}}Outer AIRGA Itn\\  $(z)$\end{tabular}}&  \textbf{\begin{tabular}[c]{@{}c@{}}Exp Pnt 1\\  $(i = 1)$\end{tabular}} &  \textbf{\begin{tabular}[c]{@{}c@{}}Exp Pnt 2\\  $(i = 2)$\end{tabular}}&  \textbf{\begin{tabular}[c]{@{}c@{}}Exp Pnt 3\\  $(i = 3)$\end{tabular}} & \textbf{$\cdots$} &  \textbf{\begin{tabular}[c]{@{}c@{}}Exp Pnt $l$\\  $(i = l)$\end{tabular}} \\ \hline
		1 & $s_{1}^{(1)}$ & $s_{2}^{(1)}$ & $s_{3}^{(1)}$ & $\cdots$ & $s_{l}^{(1)}$ \\ \hline
		2 & $s_{1}^{(2)}$ & $s_{2}^{(2)}$ & $s_{3}^{(2)}$ & $\cdots$ & $s_{l}^{2}$ \\ \hline
	\end{tabular}
\end{table}
In AIRGA, this update to the preconditioner can be done in two ways. To understand these ways, let us look at how  expansion points change in AIRGA. Table~\ref{change_in_expansion_pnt} shows changing  expansion points (labeled as Exp Pnt) for two iterations of the while loop at line 3 of Algorithm~\ref{Algo:AIRGA}. First, we can update the preconditioner when expansion points change from $s_{1}^{(1)}$ to $s_{2}^{(1)}$, $s_{2}^{(1)}$ to $s_{3}^{(1)}$, $s_{3}^{(1)}$ to $s_{4}^{(1)}$ and so on (horizontal update). Second, we can update the preconditioner when expansion points change from $s_{1}^{(1)}$ to $s_{1}^{(2)}$, $s_{2}^{(1)}$ to $s_{2}^{(2)}$, $s_{3}^{(1)}$ to $s_{3}^{(2)}$ and so on (vertical update). Since in AIRGA, vertical change in expansion points is less (which means vertically the coefficient matrices are close), we use this strategy. Thus, we are updating preconditioner from $\mathcal{K}_i^{(z-1)}  =\left(s_{i}^{(z-1)}\right)^2M + s_{i}^{(z-1)} D + K$ (which is $\mathcal{K}_{old}$) 
to $\mathcal{K}_i^{(z)}  = \left(s_{i}^{(z)}\right)^2M + s_{i}^{(z)} D + K$ (which is $\mathcal{K}_{new}$).

The AIRGA algorithm with the SPAI update preconditioner is given in Algorithm~\ref{ALGO:AIRGASPAI_Update}. Here, we only show those parts of the AIRGA algorithm that require changes. At line 10 of Algorithm~\ref{ALGO:AIRGASPAI_Update}, after solving for $Q_i^{(z)} $, ideally one would obtain the preconditioner at the $z^{\text{th}}$ AIGRA iteration as $P_i^{(z)} = Q_i^{(z)} P_i^{(z-1)}$. Since this involves a matrix-matrix multiplication, it would be expensive and defeat the purpose of using a SPAI update. Instead, we never explicitly build $P_i^{(z)}$ except at the first AIRGA iteration where we directly obtain $P_i^{(1)}$ (see line 6 of Algorithm~\ref{ALGO:AIRGASPAI_Update}). From the second AIRGA iteration onwards, we pass the all $Q_i$'s (until one reaches $P_i^{(1)}$) and $P_i^{(1)}$  to the Krylov solver so that we only require matrix-vector products (see line 11 of Algorithm~\ref{ALGO:AIRGASPAI_Update}).

%For notational ease, in the algorithm we drop the $z$ superscript. That is, $\mathcal{K}_i^{(z-1)}$ is $\mathcal{K}_{i_{old}}$ and $\mathcal{K}_i^{(z)}$ is $\mathcal{K}_i$.
%end subsection of 3
% section 4
\section{Stability Analysis of AIRGA}\label{stablity}	
An algorithm $\tilde{f}$ for computing the solution of a continuous problem $f$ on a digital computer
is said to be stable~\cite{trefethen1997numerical} if
\begin{align*}
\tilde{f}(x) = f(\tilde{x})  \ \text{for some} \  \tilde{x}\ \text{with} \ \ \frac{||\tilde{x}-x||}{||x||}= \mathcal{O}(\epsilon_{machine}),
\end{align*}              
where $\epsilon_{machine}$ is the machine precision. Here, we study the stability of AIRGA algorithm with respect to the errors introduced by iterative methods. 
\par Suppose $X^{(j)}(s_{i})$ at line 5 and 14 in AIRGA algorithm (Algorithm~\ref{Algo:AIRGA}) are computed using a direct method of solving linear system. This gives us the matrix $V$ at line 31 in Algorithm~\ref{Algo:AIRGA}. Let $f$ be the functional representation of the moment matching process that uses $V$ in AIRGA (i.e., exact AIRGA). Similarly, suppose $X^{(j)}(s_{i})$ at line 5 and 14 in Algorithm~\ref{Algo:AIRGA} are computed using an iterative method of solving linear systems. Since iterative methods are inexact, i.e., they solve the linear systems upto a certain tolerance, we denote the resulting matrix $V$ as $\tilde{V}$.      
Let $\tilde{f}$ be the functional representation of the moment matching process that uses $\tilde{V}$ in AIRGA (i.e., inexact AIRGA).  
Then, we will say that AIRGA is stable with respect to iterative solvers if%~\cite{trefethen1997numerical} 
\begin{align}\label{fst_cond}
&\tilde{f}(H(s))=f( \tilde{H}(s)) \ \ \text{for some} \ \tilde{H}(s) \ \text{with} \\ \label{secd_cond}
&\frac{\|H(s)-\tilde{H}(s) \|_{H_{2} \ \text{or} \ H_{\infty}}} {\|H(s)\|_{H_{2} \ \text{or} \ H_{\infty}}}  = \ \mathcal{O}( ||Z||),
\end{align}
where $\tilde{H}(s)$ is a perturbed original full model corresponding to the error in the linear solves for computing $\tilde{V}$ in inexact AIRGA. This perturbation is denoted by $Z$. Further, we denote $f(H(s))=\hat{H}(s)$ and $\tilde{f}(H(s))=\tilde{\hat{H}}(s).$%~\cite{ChoudharyAhuja}.  
\par In Algorithm~\ref{Algo:AIRGA}, the  linear systems at line $5$ are computed for different expansion points as
\begin{align*}
(s_{i}^{2}M+s_{i}D+K) X^{(0)}(s_{i}) = F,
\end{align*}
where $s_{i} \in \{s_{1},\ s_{2},\ \ldots,\ s_{l}\}$. As discussed earlier, we solve these linear systems inexactly  (i.e., by an iterative method). Let the residual associated with inexact linear solves for computing $X^{(0)}(s_{i})$ be $\eta_{0i}$ for $i=1,\ \ldots,\ l$ 
\begin{align}\label{inexact1}
(s_{i}^{2}M+s_{i}D+K) X^{(0)}(s_{i})&= F+\eta_{0i}.
\end{align}
%However, the similar systems can be solved exactly as 
Further, in Algorithm 1 at line 11, $\tilde{V}_{1}$ is computed as
\begin{align}\label{v_oth}
\tilde{V}_{1} =  [X^{(0)}(s_{t_{0}})/ ||X^{(0)}(s_{t_{0}})||],
\end{align}   
%where $X^{(0)}(s_{t_{0}})$ is the solution of linear systems
where $s_{t_{0}}$ is the  expansion point corresponding to the maximum moment error of the reduced system.
\par Solving the linear systems for $X^{(j)}, j= 1, \ \ldots, \ J-1$ at line 14 in Algorithm 1 inexactly yields 
\begin{align}\label{inexcat2}
(s_{i}^{2}M+s_{i}D+K) X^{(j)}(s_{i})= M\tilde{V}_{j}+\eta_{ji} \quad \text{for} \ i = 1,\ldots,l.
\end{align}
Next, $\tilde{V}_{j+1}$ is computed as   
\begin{align}\label{v_jth}
\tilde{V}_{j+1} = [X^{(j)}(s_{t_{j}})/ ||X^{(j)}(s_{t_{j}})|| ],
\end{align}
%where $X^{(j)}(s_{t_{j}})$ is the solution of linear system
where  $s_{t_{j}}$ is the expansion point corresponding to the maximum moment error of the reduced system. 

\par Finally, Galerkin projection is used to generate the reduced model (obtained by inexact AIRGA) 
\begin{align}\label{inexact_sol}
\begin{split}
\tilde{\hat{M}} = \tilde{V}^{T}M \tilde{V},\ \tilde{\hat{D}} = \tilde{V}^{T}D \tilde{V}, \ 
\tilde{\hat{K}} = \tilde{V}^{T}K \tilde{V}, \\ 
\tilde{\hat{F}}=\tilde{V}^TF, \ \tilde{\hat{C}}_p = C_p\tilde{V}, \ \text{and} \ \tilde{\hat{C}}_v = C_v\tilde{V}, 
\end{split}
\end{align}
where $\tilde{V} = [\tilde{V}_1, \ \tilde{V}_2, \ldots, \ \tilde{V}_{J}]$.
(\ref{inexact_sol}) states $\tilde f(H(s))$. Now we have to find a perturbed original model whose exact solution,\ $f(\tilde{H}(s))$,\
will give the reduced model as obtained by the inexact solution of the original full model,\ $\tilde{f}(H(s))$. That is, find $\tilde{H}(s)$ such that $\tilde{f}(H(s))=f(\tilde{H}(s))$. This would satisfy the first condition of stability~(\ref{fst_cond}).  

Assume that $\tilde H$ is given by the original matrices $M$ and $D$ and a perturbed matrix $\tilde K = K + Z$. Then, for $\tilde H$ we have 
\begin{align}\label{exact1}
(s_{i}^{2}M+s_{i}D+ (K+ Z)) X^{(0)}(s_{i})&= F \ \text{for} \ i = 1,\ldots,l.
\end{align}
%where $j=\{1,\ \ldots,\ r-1\}$
%However, similar linear systems can be solved exactly with perturbation in $K$ as
Further, assume that the linear systems can be solved exactly as
%Further, linear systems ($j=1,\ldots,r-1$) at line number 13 in Algorithm 1, can be solved exactly as  
\begin{align}\label{exact2}
(s_{i}^{2}M+s_{i}D+ (K+Z)) X^{(j)}(s_{i})=M\tilde{V}_{j} \ \text{for} \ j=1,\ \ldots,\ J-1 \ \text{and} \ i = 1,\ldots,l. 
\end{align}

Again, $\tilde{V} = [\tilde{V}_1, \ \tilde{V}_2, \ldots, \ \tilde{V}_{J}]$ where $\tilde{V}_{1}$ and $\tilde{V}_{j+1}$ for $j=1,\ 
\ldots,\ J-1$ are given by~(\ref{v_oth}) and~(\ref{v_jth}) since $X^{(0)}(s_{i})$ and $X^{(j)}(s_{i})$ for $j=1,\ \ldots,\ J-1$ 
and $i = 1,\ldots,l$ are kept same as in (\ref{inexact1}) and (\ref{inexcat2}), respectively. 

As earlier, applying Galerkin projection to the perturbed original system gives 
\begin{align}\label{exact_sol}
\begin{split}
& \hat{M} = \tilde{V}^{T}M \tilde{V},\ \hat{D} = \tilde{V}^{T}D \tilde{V}, \ \hat{K} = \tilde{V}^{T} (K+Z)\tilde{V} =  \tilde{\hat{K}} + \tilde{V}^{T} Z \tilde{V}, \\ & \hat{F} = \tilde{V}^{T}F, \ \hat{C}_{p} = C_{p}\tilde{V},\ \text{and} \ C_{v} = C_{v}\tilde{V}.
\end{split}
\end{align}
\par Our goal now is to find $Z$ such that $\hat{K}=\tilde{\hat{K}}$ (recall, that we assumed that the error can be attributed solely
to $K$; $M$ and $D$ do not change).
% because all the other reduced matrices are same~(see~(\ref{inexact_sol}) and~(\ref{exact_sol})).
Comparing~(\ref{inexact1}) with~(\ref{exact1}), and~(\ref{inexcat2}) with~(\ref{exact2}), we get
%From (\ref{inexcat2}) and (\ref{exact2}), we get
\begin{align*}
Z \ X^{(0)} (s_{t_{0}}) = - \eta_{0t_{0}} \quad \text{and} \quad
Z \ X^{(j)}(s_{t_{j}}) = -\eta_{jt_{j}} \quad \text{for} \ j = 1,\ \ldots,\ J-1,
\end{align*}
where $\eta_{0t_{0}}$ and $\eta_{jt_{j}}$ are residuals corresponding to the maximum moment error of the reduced system in inexact AIRGA.
We can rewrite $ Z $ as
\begin{align*}%\label{Orth_X}
Z \ \mathbf{X} &= -\eta,
\end{align*}
where $Z \in \mathbb{R}^{n\times n},$  $\mathbf{X}= [X^{(0)}(s_{t_0}),\ \ldots, \ X^{(J-1)}(s_{t_{(J-1)}}) ]\in \mathbb{R}^{n \times mJ},$ and $\eta = [\eta_{0t_0},\ \ldots,  \eta_{(J-1)t_{(J-1)})}]$  $\in \mathbb{R}^{n \times mJ}.$ As discussed in Section \ref{sec:1}, the upper bound for $J$ is $\lceil r/m \rceil$, and hence, $mJ < r$. Using the fact that $r \ll n$, we have $mJ < n$. Thus, we have an under-determined system of equations. 

One solution of this is 
\begin{align}\label{delta_K}
Z &=-\eta \mathbf{X}^{T}(\mathbf{X}\mathbf{X}^{T})^{-1}.
\end{align}
Multiplying both sides of (\ref{delta_K}) with $\tilde{V}$, we get
\begin{align}\label{V_oth_Z}
\tilde{V}^{T}Z \tilde{V} = - \tilde{V}^{T} \eta \mathbf{X}^{T}(\mathbf{X}\mathbf{X}^{T})^{-1}\tilde{V}.
\end{align}
For Ritz-Galerkin based iterative solvers, the solution space of linear systems is orthogonal to the residuals, i.e., $\tilde{V}_{1} \perp \eta_{0t_{0}}, \ \tilde{V}_{2} \perp \eta_{1t_{1}},\ \ldots \ ,\text{and} \ \tilde{V}_{J}  \perp  \eta_{(J-1)t_{(J-1)}}$~\cite{van2003iterative}. Hence, 
\begin{align*}
\tilde{V}^{T} \eta &= \begin{bmatrix}
\tilde{V}_{1}^{T} \\
\tilde{V}_{2}^{T} \\
\vdots &  \\
\tilde{V}_{J-1}^{T} \\
\tilde{V}_{J}^{T}
\end{bmatrix}\begin{bmatrix}
\eta_{0t_0} &
\eta_{1t_1} &
\ldots &  
\eta_{(J-1)t_{(J-1)}}   
\end{bmatrix}   \\ &=  \begin{bmatrix}  0 &  \tilde{V}^{T}_{1}\eta_{1t_{1}} & \ldots &  \tilde{V}^{T}_{1}\eta_{(J-2)t_{(J-2)}}&  \tilde{V}^{T}_{1}\eta_{(J-1)t_{(J-1)}}  \\    \tilde{V}^{T}_{2}\eta_{0t_{0}} &  0 & \ldots &  \tilde{V}^{T}_{2}\eta_{(J-2)t_{(J-2)}} & \tilde{V}^{T}_{2}\eta_{(J-1)t_{(J-1)}}   \\   \vdots & \vdots & \vdots & \vdots &\vdots\\
\tilde{V}^{T}_{J-1}\eta_{0t_{0}}  & \tilde{V}^{T}_{J-1}\eta_{1t_{1}} & \ldots & 0 &  \tilde{V}^{T}_{J-1}\eta_{(J-1)t_{(J-1)}}\\  \tilde{V}^{T}_{J}\eta_{0t_{0}} &\tilde{V}^{T}_{J}\eta_{1t_{1}} & \ldots & \tilde{V}^{T}_{J}\eta_{(J-2)t_{(J-2)}} & 0 \end{bmatrix}.
\end{align*}

Further, $\tilde{V}^{T} \eta \mathbf{X}^{T}$
\begin{align*}
	%\tilde{V}^{T} \eta \mathbf{X}^{T}\tilde{V}
 =
\begin{bmatrix}  0 &  \tilde{V}^{T}_{1}\eta_{1t_{1}} & \ldots &  \tilde{V}^{T}_{1}\eta_{(J-2)t_{(J-2)}}&  \tilde{V}^{T}_{1}\eta_{(J-1)t_{(J-1)}}  \\    \tilde{V}^{T}_{2}\eta_{0t_{0}} &  0 & \ldots &  \tilde{V}^{T}_{2}\eta_{(J-2)t_{(J-2)}} & \tilde{V}^{T}_{2}\eta_{(J-1)t_{(J-1)}}   \\   \vdots & \vdots & \vdots & \vdots &\vdots\\
\tilde{V}^{T}_{J-1}\eta_{0t_{0}}  & \tilde{V}^{T}_{J-1}\eta_{1t_{1}} & \ldots & 0 &  \tilde{V}^{T}_{J-1}\eta_{(J-1)t_{(J-1)}}\\  \tilde{V}^{T}_{J}\eta_{0t_{0}} &\tilde{V}^{T}_{J}\eta_{1t_{1}} & \ldots & \tilde{V}^{T}_{J}\eta_{(J-2)t_{(J-2)}} & 0 \end{bmatrix}
\begin{bmatrix}
X^{(0)}(s_{t_{0}})^{T} \\
X^{(1)}(s_{t_{1}})^{T} \\
\vdots   \\  X^{(J-2)}(s_{t_{J-2}})^{T} \\ 
X^{(J-1)}(s_{t_{J-1}})^{T}
\end{bmatrix}
\\ \quad \\    
=  \begin{bmatrix} 0\cdot X^{(0)}(s_{t_0})^{T} +  \tilde{V}^{T}_{1}\eta_{1t_{1}} X^{(1)}(s_{t_{1}})^{T} +  \cdots  +   \tilde{V}^{T}_{1}\eta_{(J-1)t_{(J-1)}} X^{(J-1)}(s_{t_{(J-1)}})^{T}  \\  \tilde{V}^{T}_{2}\eta_{0t_{0}} X^{(0)}(s_{t_0})^{T} +  0\cdot X^{(1)}(s_{t_1})^{T}+ \cdots \ +     \tilde{V}^{T}_{2}\eta_{(J-1)t_{(J-1)}}X^{(J-1)}(s_{t_{(J-1)}})^{T}   \\    \vdots \\  \tilde{V}_{J}^{T} \eta_{0t_{0}}X^{(0)}(s_{t_0})^{T}  +\cdots \ + \tilde{V}_{J}^{T} \eta_{(J-2)t_{(J-2)}}X^{(J-2)}(s_{t_{(J-2)}})^{T}+  0 \cdot X^{(J-1)}(s_{t_{(J-1)}})^{T}\end{bmatrix}
\end{align*}

 Since $\tilde{V}_{1} \perp \eta_{0t_{0}}, \ \tilde{V}_{2} \perp \eta_{1t_{1}},\ \ldots, \ \tilde{V}_{J}  \perp  \eta_{(J-1)t_{(J-1)}}$, and, due to (\ref{v_jth}), $\tilde{V}_{j+1}$ is just  the normalized $X^{(j)}(s_{t_j})$, we have  $X^{(0)}(s_{t_0}) \perp \eta_{0t_0}, \ X^{(1)}(s_{t_1}) \perp \eta_{1t_1},\ \ldots,$ and 
$ X^{(J-1)}(s_{t_{(J-1)}}) \perp \eta_{(J-1)t_{(J-1)}}$. Therefore from~(\ref{V_oth_Z}), we get  $\tilde{V}^{T} Z \tilde{V} =0$. 
Thus, $\hat{K} = \tilde{\hat{K}}$ or
\begin{align*}
\tilde{f}(H(s))=f(\tilde{H}(s))=\tilde{\hat{H}}(s),
\end{align*}  	
where $H(s)=(C_{p}+sC_{v})(s^{2}M+sD+K)^{-1}F,$\ $\tilde{H}(s)=(C_{p}+sC_{v})(s^{2}M+sD+ (K+ Z))^{-1}F,$ and $\tilde{\hat{H}}(s) = (\tilde{\hat{C}}_{p}+s\tilde{\hat{C}}_{v})(s^{2}\tilde{\hat{M}}+s\tilde{\hat{D}}+\tilde{\hat{K}})^{-1}\tilde{\hat{F}} = (\hat{C}_{p}+s\hat{C}_{v})(s^{2}\hat{M}+s\hat{D}+\hat{K})^{-1}\hat{F}$. Thus, we have satisfied the first condition of stability.
\par According to the second condition of stability, given in (\ref{secd_cond}), the difference between the original full model and the perturbed full model should be of the order of the perturbation~\cite{trefethen1997numerical}. This can be easily shown (Theorem 4.3 from~\cite{Beattie20122916}).
\begin{theorem}
	%\emph{(Lagrange's Theorem)}
	$ \ \text{If}  \ \  ||Z|| < \frac{1}{||\mathcal{K}(s)^{-1}||_{H_{\infty}}} \ \ \text{then} $
	\label{Second_Cond_Stab}
	\begin{align*}
	||H(s)-\tilde{H}(s)||_{H_2} \le \frac{||C(s)\mathcal{K}(s)^{-1}||_{H_{2}} ||\mathcal{K}(s)^{-1}F||_{H_{\infty}}}{1-||\mathcal{K}(s)^{-1}||_{H_{\infty}}||Z||}||Z||, 
	\end{align*} 
	where $\mathcal{K}(s) = (s^2 M + s D+ K)$ and $C(s) = (C_p + sC_v).$
\end{theorem}
Hence,
\begin{align*}
\| H(s)-\tilde{H} ( s  ) \|_{H_{2}}   = \mathcal{O}(\| Z \|).
\end{align*}
The above result holds in a relative sense too.   
This proves the stability of AIRGA. The next theorem summarizes this. 
\begin{theorem}
	%\emph{(Lagrange's Theorem)}
	\label{Tm_Stability}
	If the linear systems arising in AIRGA are solved by a Ritz-Galerkin based solver (i.e., the residual is orthogonal to the generated Krylov subspace) and   
	\begin{align*}
	||(s^2 M + s D+ K)^{-1}||_{H_{\infty}} \cdot||Z|| < 1, 	\end{align*}
	where  $Z$ given by~(\ref{delta_K}), then AIRGA is stable. 
\end{theorem}
% end Section 4
% Section 5
\section{Numerical results}\label{Sec:Numerical_exp}
%\subsection{Numerical results }
Consider a one dimensional beam model~\cite{BeattieG2005}, which is of the form (1)
\begin{align*}
& M \ddot{x}(t)+D\dot{x}(t)+Kx(t)=Fu(t), \\
& y(t)=C_px(t),
\end{align*}
where 
%$M, \ D, \  K  \in  \mathbb{R}^{n \times n} $ are mass, damping and stiffness matrices, respectively,  
$m=q=1, \ F  \in  \mathbb{R}^{n \times 1} \ \text{and} \ C_p   \in  \mathbb{R}^{1 \times n}$. The model has proportional damping, i.e., $D = \alpha M+\beta K$, where the damping coefficients $\alpha$ and $\beta$ belong to $(0,1)$~\cite{BeattieG2005}. We consider the model with two different sizes, $n = 2000$ and $n=10000.$ 

We compute a reduced order model  by the AIRGA algorithm given in Algorithm~\ref{Algo:AIRGA}. We implement AIRGA in MATLAB (2014a). We take $r_\text{max}$, i.e., the maximum dimension to which we want to reduce the system, as $30$ for model size $2000$ and $150$ for model size $10000$ based on~\cite{BeattieG2005}.  We take three expansion points that are linearly spaced between $1$ and $100$
%, and these are set as $s_1 = 1$, $s_2=500.5$ and $s_3 = 1000 $ 
based on initial data used in~\cite{BeattieG2005}. 
As discussed in Section~\ref{sec:Itr_prec}, we use iterative methods to solve the linear systems at lines 5 and 14 of Algorithm~\ref{Algo:AIRGA} instead of a direct method. Since the direct method (LU factorization) runs out of memory for large problems ($ > 50000$) even on a high configuration server (64 GB RAM), we did not pursue it further for comparison.

For the model under consideration, $M$, $D$ and $K$ matrices are symmetric positive definite. Hence, with the initial expansion points all taken as real and positive, the coefficient matrices of the linear systems to solved $s_i^2 M + s_i D+ K$ are also symmetric positive definite initially.  As discussed in Section \ref{sec:1}, after the first AIRGA iteration, the expansion points are chosen from the eigenvalues of the quadratic eigenvalue problem  $\lambda^2 \hat{M} + \lambda \hat{D} + \hat{K}$. For our example, after the first AIRGA iteration, the eigenvalues of this quadratic eigenvalue problem turn out to be complex (see Table 1.1 in \cite{Tisseur2001QuadEig} that describes why this is supposed to happen even when all matrices are symmetric positive definite). Thus, we get complex expansion points. 

Real, imaginary, or complex expansion points, each have their own merits (see Chapter 6 of \cite{grimme1997phd}). After the first AIRGA iteration, we use real parts of the complex eigenvalues as the expansion points. This is because of three reasons. First, real expansion points give good approximation for general frequency response~\cite{grimme1997phd}. Second, using real or complex expansion points has no effect on the execution of the AIRGA algorithm as well as the accuracy of the reduced system. Third, real expansion points are computationally easier to implement (the difference between $\mathcal{K}_{i}^{(z-1)}$ and $\mathcal{K}_{i}^{(z)}$ is more easily quantifiable; see Section 3.2). 

%Thus, we get complex expansion points. This implies that the coefficient matrices of the linear systems to solved after the first AIRGA iteration become non-symmetric, specifically complex-symmetric. 

For ensuring stable iterative solves in AIRGA,  from Theorem~\ref{Tm_Stability} we know that we need to use a Ritz-Galerkin based solver. Conjugate Gradient (CG) is the most popular solver based on this theory. Moreover, since CG is ideal for SPD linear systems and we obtain such linear systems here, we use CG as the underlying iterative solver. Since the unpreconditioned CG required twice as many iterations as the preconditioned CG (for almost all problem sizes), and we did not  pursue it further.

As discussed in Section~\ref{sec:Itr_prec}, preconditioning has to be employed when iterative methods fail or have very slow convergence. 
We use SPAI and SPAI update as discussed in Section~\ref{Sec:SPAI} and~\ref{Sec:SPAI_U}, respectively\footnote{For SPD linear systems, incomplete Cholesky factorization based preconditioners are quite popular. However, as discussed in Section~\ref{sec:Itr_prec}, these preconditioners are not easily parallelizable. Since for larger model sizes ($>10000$; e.g., 100,000 and so on) parallelization would be needed and incomplete Cholesky preconditioners will not be able to compete with SPAI then, we do not use them here. Moreover, SPAI preconditioners are popular even for SPD systems (see~\cite{Chow1998,George2012}).}. That is, we use Algorithm~\ref{ALGO:AIRGASPAI} with Algorithm~\ref{ALGO:SPAI} and Algorithm~\ref{ALGO:AIRGASPAI_Update} with Algorithm~\ref{ALGO:SPAI_U}. %Since for SPD linear system incomplete cholesky is the ideal preconditioner. We demonstrate its usage also.pdate 
The input to Algorithm~\ref{ALGO:SPAI} and Algorithm~\ref{ALGO:SPAI_U} is $tol$ (besides the coefficient matrices), which we take as $0.01$ based upon experience.
% Since $n_i$ is usually chosen heuristically, we use a more certain criteria. That is, residual of the minimization problem for each column less than $0.01$. This is linear cost (for each column) and we are doing such computation anyways.  
  
In Algorithms~\ref{ALGO:AIRGASPAI} and~\ref{ALGO:AIRGASPAI_Update}, at lines 1 and 2, respectively, the overall iteration (while-loop) terminates when the change in the reduced model (computed as $H_{2}$-error between the reduced models at two consecutive AIRGA iterations) is less than a certain tolerance. We take this tolerance to be $10^{-06}$ based on values in~\cite{Bonin20161}. There is one more stopping criteria in these algorithms, at lines 8 and 15, respectively.  This checks the $H_{2}$-error between two temporary reduced models. We take this tolerance to be $10^{-06}$ based on values in~\cite{Bonin20161}. Since this is an adaptive algorithm, the optimal size of the reduced model is determined by the algorithm itself, and is denoted by \textit{r}.
\begin{table}[]
	%[!tbp]
	%\normalsize
	\centering
	%\captionsetup{font=scriptsize}
	\caption{ CG iterations and computation time for model size 2000}
	\label{GMRES_SPAI_500}
	\def\arraystretch{1.1}
	\begin{tabular}{|c|c|c|l|c|}
		\hline
		\multirow{2}{*}{\textbf{AIRGA  Iteration\#}} & \multicolumn{2}{c|}{\textbf{\begin{tabular}[c]{@{}c@{}}CG using SPAI \end{tabular}}} & \multicolumn{2}{c|}{\textbf{\begin{tabular}[c]{@{}c@{}}CG using  SPAI Update \\ \end{tabular}}} \\ \cline{2-5} 
		& \textbf{Iter}                                       & \textbf{Time (secs)}    & \textbf{Iter}    & \multicolumn{1}{l|}{\textbf{Time (secs)}}    \\ \hline	\multirow{3}{*}{1}  & 4   & 0.07   & 4      & 0.07         \\ \cline{2-5} & 4       & 0.07    & 4    & 0.07         \\ \cline{2-5} & 4  & 0.07    & 4        & 0.07        \\ \hline \multirow{3}{*}{2}  & 3     & 0.06     & 3     & 0.06   \\ \cline{2-5} 	& 3        & 0.06   & 3    & 0.06   \\ \cline{2-5} & 4          & 0.07     & 4    & 0.07     \\ \hline
		\multirow{3}{*}{3}  & 4     & 0.07     & 4     & 0.11  \\ \cline{2-5} 	& 4        & 0.07   & 4     & 0.11   \\ \cline{2-5} & 4          & 0.07     & 4   & 0.11     \\ \hline
		\multirow{3}{*}{4}  & 4     & 0.07     & 4     & 0.18  \\ \cline{2-5} 	& 4        & 0.07   & 4     & 0.18   \\ \cline{2-5} & 4          & 0.07     & 4    & 0.18     \\ \hline
	\end{tabular}
\end{table}

\begin{table}[]
	%[!htbp]
	\centering
	%\normalsize
	\caption{ CG iterations and computation time for model size 10000}
	\label{GMRES_SPAI_1000}
	\def\arraystretch{1.1}
	\begin{tabular}{|c|c|c|l|c|}
		\hline
		\multirow{2}{*}{\textbf{AIRGA Iteration\#}} & \multicolumn{2}{c|}{\textbf{\begin{tabular}[c]{@{}c@{}}CG using SPAI \end{tabular}}} & \multicolumn{2}{c|}{\textbf{\begin{tabular}[c]{@{}c@{}}CG using SPAI update \\ \end{tabular}}} \\ \cline{2-5} 
		& \textbf{Iter} & \textbf{Time (secs)} & \textbf{Iter} & \multicolumn{1}{l|}{\textbf{Time (secs)}} \\ \hline
		\multirow{3}{*}{1} & 3 & 1.4 & 3 & 1.4 \\ \cline{2-5} 
		& 4 & 1.4 & 4 & 1.5 \\ \cline{2-5} 
		& 4 & 1.5 & 4 & 1.5 \\ \hline
		\multirow{3}{*}{2} & 4 & 1.5 & 4 & 1.5 \\ \cline{2-5} 
		& 3 & 1.4 & 3 & 1.4 \\ \cline{2-5} 
		& 3 & 1.4 & 3 & 1.4 \\ \hline
		\multirow{3}{*}{3} & 4 & 1.5 & 4 & 2.4 \\ \cline{2-5} 
		& 3 & 1.4 & 3 & 2.3 \\ \cline{2-5} 
		& 4 & 1.5 & 4 & 2.4 \\ \hline
		\multirow{3}{*}{4} & 4 & 1.4 & 4 & 4.1 \\ \cline{2-5} 
		& 3 & 1.4 & 3 & 3.7 \\ \cline{2-5} 
		& 4 & 1.4 & 4 & 4.1 \\ \hline
	\end{tabular}
\end{table}
\begin{table}[]
	%[!tp]
	%\normalsize
	\centering
	\caption{ SPAI and SPAI update computation time for model size 2000}
	\label{SPAI_SPAI_Update_time_2K}
	\def\arraystretch{1.1}
	\begin{tabular}{|c|c|c|}%{p{3.5 cm} p{8 cm} p{5 cm}}
		\hline
		\textbf{AIRGA  Iteration\#} & \textbf{\begin{tabular}[c]{@{}c@{}} SPAI (secs)\end{tabular}} & \textbf{\begin{tabular}[c]{@{}c@{}}  SPAI with update (secs)\end{tabular}} \\ \hline
		\multirow{3}{*}{1} & 10.0 & 10.0 \\ \cline{2-3} 
		& 35.0 & 35.0 \\ \cline{2-3} 
		& 36.0 & 36.0\\ \hline
		\multirow{3}{*}{2}  & 10.0 & 10.0 \\ \cline{2-3} 
		& 11.0 & 11.0\\ \cline{2-3} 
		& 10.0 & 10.0\\ \hline
		\multirow{3}{*}{\textbf{3}} & \textbf{10.0} &\textbf{ 0.8} \\ \cline{2-3} 
		& \textbf{11.0} & \textbf{6.0} \\ \cline{2-3} 
		& \textbf{11.0} & \textbf{1.1} \\ \hline
		\multirow{3}{*}{\textbf{4}}  & \textbf{11.0} & \textbf{0.6} \\ \cline{2-3} 
		& \textbf{10.0} & \textbf{4.0 }\\ \cline{2-3} 
		& \textbf{11.0 }&\textbf{ 0.8}\\ \hline
	\end{tabular}
\end{table}

\begin{table}
	%[!bp]
	\centering
	%\normalsize
	\caption{SPAI and SPAI update computation time for model size 10000}
	\label{SPAI_Update_1000}
	\def\arraystretch{1.1}
	\begin{tabular}{|c|c|c|}
		\hline
		\textbf{AIRGA Iteration\#} & \textbf{\begin{tabular}[c]{@{}c@{}} SPAI (secs)\end{tabular}} & \textbf{\begin{tabular}[c]{@{}c@{}}  SPAI  with  Update (secs)\end{tabular}} \\ \hline
		\multirow{3}{*}{1}  & 282.0 & 282.0 \\ \cline{2-3} 
		& 476.0 & 476.0 \\ \cline{2-3} 
		& 530.0 & 530.0 \\ \hline
		\multirow{3}{*}{2}  & 177.0 & 177.0 \\ \cline{2-3} 
		& 276.0 & 276.0 \\ \cline{2-3} 
		& 170.0 & 170.0\\ \hline
		\multirow{3}{*}{\textbf{3}}  & \textbf{172.0} & \textbf{28.0}\\ \cline{2-3} 
		& \textbf{272.0} & \textbf{146.0} \\ \cline{2-3} 
		& \textbf{280.0} & \textbf{37.0}\\ \hline
		\multirow{3}{*}{4}  &\textbf{ 179.0} & \textbf{27.0} \\ \cline{2-3} 
		& \textbf{300.0} & \textbf{65.0}\\ \cline{2-3} 
		&\textbf{ 179.0} & \textbf{22.0}\\ \hline
	\end{tabular}
\end{table}

\begin{table}[]
	\centering
	\caption{SPAI and SPAI update  analysis for model size 2000}
	\subfloat[]{
	\label{norm_anyl_2K_ep_1}
	\begin{tabular}{|c|c|c|c|}
		\hline
		\multicolumn{1}{|c|}{\multirow{2}{*}{\textbf{AIRGA Iteration ($z$)}}} & \multicolumn{3}{c|}{\textbf{Expansion Points 1}} \\ \cline{2-4}
		\multicolumn{1}{|c|}{} & \textit{\textbf{Value}} & \textit{\textbf{$\left|\left|I-\mathcal{K}_1^{(z)}\right|\right|_F$}} & \textit{\textbf{$\left|\left|\mathcal{K}_{1}^{(z-1)}- \mathcal{K}_{1}^{(z)}\right|\right|_F$}} \\[5 pt] \hline
%		1 & 3.14 & $   1.5e+03$ & \textit{NA} \\ \hline
%		2 & 0.2688 & $110.59$ & \textit{NA} \\ \hline
		3 & 0.2681 & $  110.5137
		$ & $    0.1124
		$ \\ \hline
		4 & 0.2682 & $110.5170$ & $0.0046$ \\ \hline
	\end{tabular}}
	\vspace{0.2cm}
	%\end{table}
	%\begin{table}[]
	\centering
	\subfloat[]{
	\label{norm_anyl_2K_ep_2}
	\begin{tabular}{|c|c|c|c|}
		\hline
		\multicolumn{1}{|c|}{\multirow{2}{*}{\textbf{AIRGA Iteration ($z$)}}} & \multicolumn{3}{c|}{\textbf{Expansion Point 2}} \\ \cline{2-4} 
		\multicolumn{1}{|c|}{} & \textit{\textbf{Value}} & \textit{\textbf{$\left|\left|I- \mathcal{K}_{2}^{(z)}\right|\right|_F$}} & \textit{\textbf{$\left|\left|\mathcal{K}_{2}^{(z-1)}- \mathcal{K}_{2}^{(z)}\right|\right|_F$}} \\[5 pt] \hline
%		1 & 158.65 & $3.0e+06$ & \textit{NA} \\ \hline
%		2 & $2.0640$ & $741.10$ & \textit{NA} \\ \hline
		3 & $1.9948$ & $701.4626$ & $40.2943$ \\ \hline
		4 & $1.8491$ & $  631.3391
		 $ & $9.6520$ \\ \hline
	\end{tabular}}
	\vspace{0.2cm}
	%\end{table}
	%\begin{table}[]
	\centering
	\subfloat[]{
	\label{norm_anyl_2K_ep_3}
	\begin{tabular}{|c|c|c|c|}
		\hline
		\multicolumn{1}{|c|}{\multirow{2}{*}{\textbf{AIRGA Iteration ($z$)}}} & \multicolumn{3}{c|}{\textbf{Expansion Point 3}} \\ \cline{2-4} 
		\multicolumn{1}{|c|}{} & \textit{\textbf{Value}} & \textit{\textbf{$\left|\left|I-\mathcal{K}_{3}^{(z)}\right|\right|_F$}} & \textit{\textbf{$\left|\left|\mathcal{K}_{3}^{(z-1)}- \mathcal{K}_{3}^{(z)}\right|\right|_F$}} \\[5 pt]  \hline
%		1 & 314.2 & $1.2e+07$ & \textit{NA} \\ \hline
%		2 & 0.2740 & $111.2016$ & \textit{NA} \\ \hline
		3 & 0.2700 & $ 110.6989$ & $0.0398$ \\ \hline
		4 & 0.2699 & 110.7225 & $  0.0068$ \\ \hline
	\end{tabular}}
\end{table}

\begin{table}[]
	\centering
	\caption{SPAI and SPAI update  analysis for model size 10000}
	\subfloat[]{
	\label{norm_anyl_10K_ep_1}
	\begin{tabular}{|c|c|c|c|}
		\hline
		\multicolumn{1}{|c|}{\multirow{2}{*}{\textbf{AIRGA Iteration ($z$)}}} & \multicolumn{3}{c|}{\textbf{Expansion Points 1}} \\ \cline{2-4} 
		\multicolumn{1}{|c|}{} & \textit{\textbf{Value}} & \textit{\textbf{$\left|\left|I-\mathcal{K}_1^{(z)}\right|\right|_F$}} & \textit{\textbf{$\left|\left|\mathcal{K}_{1}^{(z-1)}- \mathcal{K}_{1}^{(z)}\right|\right|_F$}} \\[5 pt]  \hline
%		1 & 3.14 & $3.4e+03$ & \textit{NA} \\ \hline
%		2 & 0.269 & $247.48$ & \textit{NA} \\ \hline
		3 & 0.2681 & $247.18$ & $0.4279$ \\ \hline
		4 & 0.2682 & $247.17$ & $0.0066$ \\ \hline
	\end{tabular}}
	\vspace{0.2cm}
	%\end{table}
	%\begin{table}[]
	\centering
	\subfloat[]{
	\label{norm_anyl_10K_ep_2}
	\begin{tabular}{|c|c|c|c|}
		\hline
		\multicolumn{1}{|c|}{\multirow{2}{*}{\textbf{AIRGA Iteration ($z$)}}} & \multicolumn{3}{c|}{\textbf{Expansion Point 2}} \\ \cline{2-4} 
		\multicolumn{1}{|c|}{} & \textit{\textbf{Value}} & \textit{\textbf{$\left|\left|I- \mathcal{K}_{2}^{(z)}\right|\right|_F$}} & \textit{\textbf{$\left|\left|\mathcal{K}_{2}^{(z-1)}- \mathcal{K}_{2}^{(z)}\right|\right|_F$}} \\[5 pt]  \hline
%		1 & 158.65 & $6.8e+06$ & \textit{NA} \\ \hline
%		2 & $1.87$ & $1.4e+03$ & \textit{NA} \\ \hline
		3 & $1.5037$ & $1.0e+03$ & $416.5$ \\ \hline
		4 & $1.9415$ & $ 1.5e+03 $ & $116.8$ \\ \hline
	\end{tabular}}
	\vspace{0.2cm}
	%\end{table}
	%\begin{table}[]
	\centering
	\subfloat[]{
	\label{norm_anyl_10K_ep_3}
	\begin{tabular}{|c|c|c|c|}
		\hline
		\multicolumn{1}{|c|}{\multirow{2}{*}{\textbf{AIRGA Iteration ($z$)}}} & \multicolumn{3}{c|}{\textbf{Expansion Point 3}} \\ \cline{2-4} 
		\multicolumn{1}{|c|}{} & \textit{\textbf{Value}} & \textit{\textbf{$\left|\left|I-\mathcal{K}_{3}^{(z)}\right|\right|_F$}} & \textit{\textbf{$\left|\left|\mathcal{K}_{3}^{(z-1)}- \mathcal{K}_{3}^{(z)}\right|\right|_F$}} \\[5 pt]  \hline
%		1 & 314.2 & $2.6e+07$ & \textit{NA} \\ \hline
%		2 & 0.280 & $250.4$ & \textit{NA} \\ \hline
		3 & 0.2698 & $247.6$ & $3.9324$ \\ \hline
		4 & 0.2696 & $247.5$ & $0.0602$ \\ \hline
	\end{tabular}}
\end{table}

\begin{table}
	%[!tbp]
	\centering
	%\normalsize
	\caption{Total computation time of CG with  preconditioners}
	\label{Total_GMRES_pre}
	\def\arraystretch{1.1}
	\begin{tabular}{|c|c|c|}
		\hline
		\textbf{Size} & \textbf{\begin{tabular}[c]{@{}c@{}} CG and SPAI  (Min)\end{tabular}} & \textbf{\begin{tabular}[c]{@{}c@{}}CG and SPAI update  (Min)\end{tabular}} \\ \hline
		2000  & 3.1                                                                               & 2.1                                                                      \\ \hline
		10000 & 61.2                                                                            & 
		40.7                                                                   \\ \hline
	\end{tabular}
\end{table}

\begin{table}
	%[!]
	\centering
	\def\arraystretch{1.1}
	%\normalsize
	\caption{Accuracy of reduced system}
	\label{Accuracy}
	\begin{tabular}{|l|l|c|l|l|}
		\hline
		\textbf{Problem} & \textbf{n} & \textbf{Method} & \textbf{Error} & \textbf{r} \\ \hline
		\multirow{6}{*}{\textbf{1-D Beam Model}} & \multirow{3}{*}{2000} & LU & $1.2 e-06 $ & 27\\ \cline{3-5} 
		&  & \begin{tabular}[c]{@{}c@{}}CG with SPAI\end{tabular} & $3.0 e-06 $& 27 \\ \cline{3-5} 
		&  & \begin{tabular}[c]{@{}c@{}}CG with SPAI update\end{tabular} & $ 2.2 e-06 $& 27 \\ \cline{2-5} 
		& \multirow{3}{*}{10000} & LU & $1.1e-06 $& 128 \\ \cline{3-5} 
		&  & \begin{tabular}[c]{@{}c@{}}CG with SPAI\end{tabular} & $3.6 e-06$ & 128 \\ \cline{3-5} 
		&  & \begin{tabular}[c]{@{}c@{}}CG with SPAI update\end{tabular} & $2.7 e-06$ & 128 \\ \hline
	\end{tabular}
\end{table}
As discussed in Section~\ref{Sec:SPAI_U}, SPAI update in AIRGA is done vertically (see Algorithm~\ref{ALGO:AIRGASPAI_Update} and Table~\ref{change_in_expansion_pnt}). Ideally, this update should  be done from AIRGA iteration 2 (since at iteration 1, there is no preceding set of expansion points). However, for this problem change in expansion points from iteration 1 to 2 is fairly large (implying the corresponding coefficient matrices are far). Hence, we implement update from AIRGA iteration 3 onwards.     

Table~\ref{GMRES_SPAI_500} gives the CG iteration count and time when using  basic SPAI (that is without SPAI update) and SPAI with update for model size 2000. Table~\ref{GMRES_SPAI_1000} gives the same data for model size 10000. From Tables~\ref{GMRES_SPAI_500} and~\ref{GMRES_SPAI_1000} it is observed that the CG computation time for both variants of preconditioners  is  almost same. The computation time for CG using SPAI update is slightly higher than CG using basic SPAI from AIRGA iteration 3 onwards (for both model sizes). This is because from AIRGA iteration 3 onwards, we apply SPAI update and hence, the number of matrix-vector products increase (see lines 11 and 17 of Algorithm \ref{ALGO:AIRGASPAI_Update}). However, this slight increase in the time for SPAI with update fades when one takes the preconditioner computation time into account (see the discussion below).
%fairy small.  

Table~\ref{SPAI_SPAI_Update_time_2K} gives the time for computing the basic SPAI preconditioner  and SPAI with update preconditioner  for model size 2000. Table~\ref{SPAI_Update_1000} gives the same data for model size 10000. From Tables~\ref{SPAI_SPAI_Update_time_2K} and~\ref{SPAI_Update_1000} it is observed that considerable amount of time is saved  by using SPAI with updates. As discussed in earlier paragraphs, since SPAI update is being done only from AIRGA iteration 3, saving in time is observed from this step onwards only. 

We also analyze why SPAI update takes less time. As discussed in Section~\ref{Sec:SPAI_U}, SPAI update is useful when $\left|\left|I-\mathcal{K}_{i}^{(z)}\right|\right|_F$ is large and $\left|\left|\mathcal{K}_{i}^{(z-1)} - \mathcal{K}_{i}^{(z)}\right|\right|_F$ is small. This data for three expansion points for model size 2000 is given in Tables 6(A), 6(B) and 6(C). Similar data for model size 10000 is given in Tables 7(A), 7(B) and 7(C). For model size 2000, in Table~\ref{SPAI_SPAI_Update_time_2K} we see that at AIRGA iterations 3 and 4 computation time for SPAI update is almost one fourth of the SPAI time. This is because $\left|\left|\mathcal{K}_{i}^{(z-1)} - \mathcal{K}_{i}^{(z)}\right|\right|_F$ is very small as compared to $\left|\left|I- \mathcal{K}_{i}^{(z)}\right|\right|_F$ (see Tables 6(A), 6(B) and 6(C)). Similar pattern is observed for model size 10000 (see Table~\ref{SPAI_Update_1000} and Tables 7(A), 7(B) and 7(C)).

\par Table~\ref{Total_GMRES_pre} shows the total computation time for iterative solves (i.e., CG time plus the preconditioner time) when using basic SPAI preconditioner  and when using SPAI with update preconditioner  for model sizes 2000 and 10000. 
%We can notice from this table that iterative solves with update take about ? of time as needed for basic iterative solves (or without SPAI update).
We can notice from this table that computation time of iterative solves with SPAI update is on an average $\frac{2}{3}$-rd of the computation time of iterative solves with SPAI. This saving is larger for model size 10000   (we go from around 1 hour to 40 Minutes). Hence, larger the problem more the saving.

Table~\ref{Accuracy} lists the relative error between the original model and the reduced model as well as the size to which the model is reduced ($r$) when using LU factorization (direct method), CG with SPAI, and CG with SPAI updates. Since the error and $r$ values for all the above three cases are  almost the same, we can conclude that by using iterative solves in AIRGA, the quality of reduced system is not compromised. 

\section{Conclusion and Future work}\label{sec:conl_fut} 
We discussed the application of preconditioned iterative methods for solving large linear systems arising in AIRGA. The SPAI preconditioner works well here and SPAI update (where we reuse the preconditioner) leads to substantial savings. This is demonstrated by experiments on two different sizes of one dimensional beam model. We also presented conditions under which the AIRGA algorithm is stable with respect to the errors introduced by iterative methods. 

Future work includes applying preconditioned iterative methods in other model reduction algorithms for second order dynamical systems~(besides AIRGA). For example, Alternate Direction Implicit (ADI) methods for model reduction of second order linear dynamical systems~\cite{Uddin2012}. Based upon our studies on AIRGA and ADI based methods, we also plan to propose a class of preconditioners that would work for most  model reduction algorithms for second order linear dynamical systems.
\section*{Acknowledgement}
We would like to thank Prof. Eric de Sturler (at Department of Mathematics, Virginia Tech, Blacksburg, VA, USA) for stimulated discussion regarding stability of AIRGA.
%\section*{References}

\bibliographystyle{amsplain}
\bibliography{Mybib}

\end{document}